\documentclass[a4paper,10pt]{article}

\usepackage{amsmath, amsfonts, amssymb, latexsym, mathrsfs, graphicx, nopageno}

\usepackage[all]{xy}
\newdir{ ->}{{}*!/-5ex/@{->}}

\newcommand{\mmap}[3]{\ar@/#2/[#1]|*=0{\rotatebox{#3}{$\scriptsize |$}}} %\mmap{Direction}{CurvatureNoneIfEmpty}{Rotation90IsInvariant}%
\newcommand{\picar}[3]{\ar@{}[#1]|*{\rotatebox{#2}{$\scriptsize #3$}}} %\picar{Direction}{Rotation90IsInvariant}{Symbol}%

\newcommand{\ba}{{\cal{A}}}

\newcommand{\bh}{{\cal{H}}}
\newcommand{\oop}{\operatorname{op}}
\newcommand{\op}{^{\oop}}

\newcommand{\ev}{\mathbb{V}}

\begin{document}

\begin{center}
{\underline{{\Large Enriched Herds And Finite Quantum Groupoids}}}

\begin{center}
{{\small Thomas Booker and Brian J. Day}}

{{\small April 13, 2010.}}
\end{center}
\end{center}

\noindent \emph{\small \underline{Abstract:} We describe a finite quantum groupoid associated to any finite $Vect_{k}$-enriched herd.}
\newline

\begin{center}
\begin{tabular}{c}\hline
\phantom{a}\phantom{a}\phantom{a}\phantom{a}\phantom{a}\phantom{a}\phantom{a}\phantom{a}\phantom{a}\phantom{a}\phantom{a}\phantom{a}\phantom{a}\phantom{a}\phantom{a}\phantom{a}
\end{tabular}
\end{center}

Let $\ev = Vect_{k}$ (finite dimensional $k$-vector spaces) and let $\ba$ be a finite representable $k$-linear flock (i.e. a finite $\ev$-herd): thus the objects of $\ba$ are finite in number, $\ba ( a, b)$ finite, $\ba(a,b)^{\ast} \cong \ba(b,a)$ chosen
(we refer the reader to [1] for the definition of a $k$-linear flock).

\noindent Then the category
$$[\ba\op\otimes\ba , Vect_{k}] $$

\noindent of finite bimodules over $\ba$ is \emph{$\ast$-autonomous} (as shown in [3]).

\noindent If also
$$\ba(b, \tau(a,c,d))^{\ast} \cong \ba(a,\tau (b,d,c))$$

\noindent then we obtain a functor (i.e. antipode)
\begin{eqnarray*}S: & \bh\op \longrightarrow \bh \\ & (x,y) \mapsto (y,x)\end{eqnarray*}

\noindent where the objects of $\bh$ are the pairs $(a,b)\in \ba\op\otimes\ba$ and the hom of $\bh$ is

$$\bh((a,b),(c,d))=\ba(\tau(b,a,c),d)\,\, .$$

\noindent Moreover, the convolution $[\bh,Vect_{k}]$ with respect to the promonoidal structure on $\bh$ given by

\begin{equation}p((a,b),(c,d),(u,v)) = \ba(\tau(d,c,b),v)\otimes\ba(u,a)\end{equation}
$$j(u,v) = \ba(u,v)$$

\noindent is $\ast$-autonomous for the above $S$.

There is a prospective Fourier transformation functor $\hat{K}$ on $[\bh,Vect_{k}]$ given by

$$\hat{K}(f)(a,b) = \int\limits^{(x,y)\in\bh} f(x,y)\otimes\ba(\tau(y,x,a),b)$$

\noindent where

$$\xymatrix{[\bh,Vect_{k}] \ar@<+1ex>[rr]^{\hat{K}} \ar@<-1ex>@{<-}[rr]_{\check{K}} && [\ba\op\otimes\ba,Vect_{k}]}$$

\noindent with $\hat{K}\dashv\check{K}$. 
The cocontinuous functor $\hat{K}$ is certainly multiplicative and preserves the induced duality on $[\bh,Vect_{f}]$ defined  by $f^{\ast}(x,y) = f(y,x)^{\ast}$.
Further, $\hat{K}$ is conservative (as required in [2]).

\pagebreak

In fact the Fourier functor $\hat{K}$ is here nothing other than the process of ``restriction along $h$'' where
$$h:\ba\op\otimes\ba\to\bh$$
is the canonical Kleisli functor mapping $(a,b)\in\ba\op\otimes\ba$ to $(a,b)\in\bh$, which is surjective on objects (see [1] for the details of this construction).
The functor $\hat{K}=[h,1]$ also preserves both the left and right internal homs of the convolution structure on $[\bh , Vect_{k}]$, this is a straight forward calculation.
So $[\bh , Vect_{k}]$ is a $\ast$-autonomous monoidal biclosed category because $[\ba\op\otimes\ba, Vect_{k}]$ is such (and $\hat{K}=[h,1]$ is conservative).
The $\ast$-autonomy of $[\ba\op\otimes\ba, Vect_{k}]$ also implies that the domain $\bh$ is itself $\ast$-autonomous  as a promonoidal category  over $Vect_{k}$, that is, the promultiplication $p$ on $\bh$ mentioned above (1) satisfies the cyclic relations given in [2] and [3] for example.
\newline

Thus $\hat{K}=[h,1]$ is a quantum groupoid in the sense of [3] (see also [5]), and there is undoubtedly also a connection with the article by T. Brzezi\'nski and J. Vercruysse: ``Bimodule herds'', J. Alg. 321 (2009), 2670-2704.
Indeed, the Kleisli functor (or ``algebroid'')
$$h:\ba\op\otimes\ba\to\bh$$
may be related to the Ehresmann Groupoid, even in the more general case where $\ba$ is not necessarily finite as a $\ev$-herd.

\begin{center}
\begin{tabular}{c}\hline
\phantom{a}\phantom{a}\phantom{a}\phantom{a}\phantom{a}\phantom{a}\phantom{a}\phantom{a}\phantom{a}\phantom{a}\phantom{a}\phantom{a}\phantom{a}\phantom{a}\phantom{a}\phantom{a}
\end{tabular}
\end{center}

\begin{center}
\small{Mathematics Dept., Faculty of Science, Macquarie University, NSW 2109, Australia.}
\newline

Replies are welcome through Tom Booker (thomas.booker@students.mq.edu.au).
\end{center}

\end{document}